\documentclass[twocolumn]{autart}
\let\theoremstyle\relax
\usepackage{amsmath,amssymb,paralist,comment,array,cite,enumitem,graphicx,amsthm,xcolor,soul}
\usepackage{mathtools}
\newtheorem{theorem}{Theorem}
\newtheorem{lemma}{Lemma}

\newtheorem{remark}{Remark}
\usepackage[numbers]{natbib}
\usepackage{soul}
\usepackage{algorithm}
\usepackage{algpseudocode}
\begin{document}
\begin{frontmatter}

\title{Time-Varying Parameters in Sequential Decision Making Problems} 
\thanks{This work was supported by DOE award DE-EE0009125, and Dynamic Research Enterprise for Multidisciplinary Engineering Sciences (DREMES) - collaboration between Zhejiang University and the University of Illinois at Urbana-Champaign. This paper was not presented at any IFAC 
meeting. Corresponding author Amber Srivastava.}
\author[Sriva]{A Srivastava}\ead{asrivastava@control.ee.ethz.ch},    
\author[Salapaka]{S M Salapaka}\ead{salapaka@illinois.edu},               
\address[Sriva]{Automatic Control Laboratory, Swiss Federal Institute of Technology (ETH Zurich), Physicstrasse 3, 8092 Zurich, Switzerland}  
\address[Salapaka]{Coordinated Science Laboratory, and Mechanical Science and Engineering, University of Illinois Urbana-Champaign, USA}
          
\begin{abstract}                          
In this paper we address the class of Sequential Decision Making (SDM) problems that are characterized by time-varying {\em parameters}. These parameter dynamics are either {\em pre-specified} or {\em manipulable}. At any given time instant the decision policy --- that governs the sequential decisions --- along with all the parameter values determines the cumulative cost incurred by the underlying SDM. Thus, the objective is to determine the {\em manipulable} parameter dynamics as well as the {\em time-varying} decision policy such that the associated cost gets minimized at each time instant. To this end we develop a control-theoretic framework to design the unknown parameter dynamics such that it locates and tracks the optimal values of the parameters, and simultaneously determines the time-varying optimal sequential decision policy. Our methodology builds upon a Maximum Entropy Principle (MEP) based framework that addresses the static parameterized SDMs. More precisely, we utilize the resulting smooth approximation (from the above framework) of the cumulative cost as a control Lyapunov function. We show that under the resulting control law the parameters asymptotically track the local optimal, the proposed control law is Lipschitz continuous and bounded, as well as ensure that the decision policy of the SDM is optimal for a given set of parameter values. The simulations demonstrate the efficacy of our proposed methodology.
\end{abstract}
\end{frontmatter}

\section{Introduction}\label{sec: Intro}
Sequential Decision Making (SDM) problems are ubiquitous in engineering. Popular application areas such as job-shop scheduling, vehicle routing, sensor networks, and autonomous robotics involve SDM problems \cite{littman1996algorithms}. These problems are typically characterized by discrete-time dynamic control systems, which describe how one-step evolution of a {\em state} of the system depends on the control {\em action}; and a cost function that specifies the cost incurred in this one-step evolution of the state. The related objective is to determine a {\em decision policy} (a sequence of actions) that minimizes the cumulative cost incurred over finite or infinite horizon. The SDM problems are modelled using several different mathematical frameworks --- such as optimal control, dynamic programming \cite{liberzon2011calculus}, Markov decision processes (MDPs) \cite{bertsekas1996neuro}, probabilistic automaton \cite{stoelinga2002introduction}, and model predictive control (MPC) \cite{camacho2013model} --- where prior literature provides extensive solution methodologies such as Pontryagin's maximum principle \cite{liberzon2011calculus}, value and policy iteration, linear programming, and reinforcement learning \cite{bertsekas1996neuro,bertsekas2011dynamic}.

\begin{figure}
\centering
\includegraphics[width=\columnwidth]{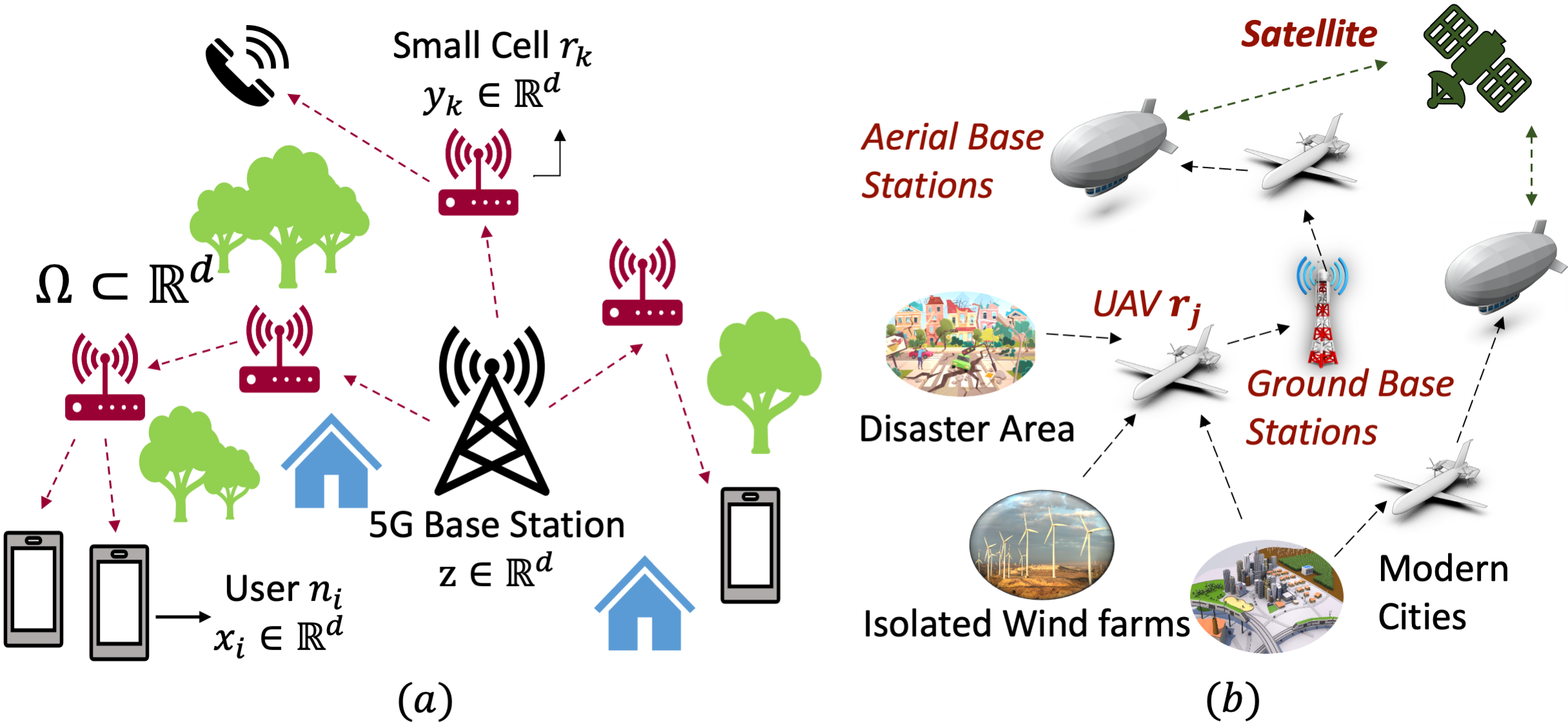}
\caption{(a) 5G small cell network - static para-SDM. (b) Multi-UAV network - dynamic para-SDM.}\label{fig: 5Gsmall_cell}
\end{figure}

Scenarios such as self organizing networks \cite{aguilar2015location}, 5G small cell networks \cite{siddique2015wireless}, supply chain, UAV communcation systems \cite{shakeri2019design}, and last mile delivery problems \cite{9096570} pose a {\em parameterized} Sequential Decision Making ({\em para}-SDM) problem. Here the main difference is that the states and actions, and the cost function themselves depend on external parameters. Some of these parameters may themselves be manipulable, and form a part of the decision variables for the underlying optimization problem. The objective in these para-SDM problems are to simultaneously (a) determine the optimal decision policy governing the sequential decision making, as well as (b) ascertain the unknown (or manipulable) parameters in the problem such that the associated cumulative cost gets minimized. For instance, consider the 5G-small cell network illustrated in the Figure \ref{fig: 5Gsmall_cell}(a). Here, the users $\{n_i\}$ distributed at the locations $\{x_i\in \mathbb{R}^d\}$ are required to communicate back and forth with the base station $\delta$ located at $z\in\mathbb{R}^d$. The objective is to simultaneously (a) overlay a network of small cells $\{r_j\}$ on the existing network of users and base station, and (b) determine the route (possibly multi-hop) between the users and the base station via the network of small cells such that the total communication cost (for instance, total network delay) gets minimized. Here, the unknown small cell locations $\{y_j\}$ and the shortest routes, respectively, constitute the unknown parameters and the decision policy for the underlying para-SDM. In particular, at each stage of the route the decision policy determines the next state --- the next small cell --- where each such state $r_j$ is parameterized by its location $y_j$. In the context of our current work, we refer to such problems as {\em static} para-SDM. These problems come with a lot of inherent complexities. For instance, the latter objective of determining the parameters (small cell locations) is akin to facility location problem that is shown to be NP-hard \cite{mahajan2009planar} with a non-convex cost surface riddled with multiple poor local minima. Further, due to the additional state and action parameters it is difficult to model para-SDMs directly by using the existing frameworks \cite{liberzon2011calculus,bertsekas1996neuro,stoelinga2002introduction,camacho2013model} that model SDMs. We have addressed the static para-SDMs in \cite{9517030}.

In this article, we address the {\em dynamic} para-SDM problems. Here, the parameters have an associated dynamics. The dynamics of some of these parameters maybe manipulable, i.e., the parameter dynamics is represented by a control system $\dot \Upsilon=f(\Upsilon,u,t)$, where the feedback law $u(\Upsilon, t)$ needs to be designed (along with determining the decision policy) such that the cost objective representing the dynamic para-SDM is minimized at each time instant. For instance, in the context of 5G-small cell networks the user nodes $\{n_i\}$ could be mobile, i.e. their respective spatial locations $\{x_i\in\mathbb{R}^d\}$ change with time. As a result, the initial optimal routes (governed by the decision policy) and the small cell locations $\{y_j \in \mathbb{R}^d\}$ --- that minimize the total network cost at the time $t=0$ --- will no longer be optimal for the future time instants $t>0$. In other words, the communication routes and the small cell locations need to evolve with time so that the network cost gets minimized at {\em each} time instant $t$. This can be further understood in terms of a similar multi-UAV (unmanned aerial vehicle) communication network (see Figure \ref{fig: 5Gsmall_cell})(b) \cite{shakeri2019design}. As illustrated in the Figure, the network utilizes a multi-UAV system to effectively route the communication packets to (and from) the aerial base station $\delta$, and at the same time provide an appropriate coverage for its user nodes. Since the network also comprises of mobile users, the UAVs need to evolve dynamically so as to facilitate appropriate coverage and routing at all times in the network. Thus, the objective is to determine the dynamics governing the UAV locations as well as the time-varying communication routes in the wireless network that minimizes the communication cost at each time instant. The dynamic para-SDMs are a difficult class of problems as on one hand they inherit all the existing complexities of static para-SDM problems (as stated above), and on the other hand they require to determine time-varying solutions.

One of the main contributions of our earlier work \cite{9517030} on para-SDM is to view them as combinatorial optimization problems. This is owing to the combinatorially large number of possible sequences of states and actions (also, referred to as {\em paths}) in the SDM. This viewpoint enables the use of Maximum Entropy Principle (MEP) --- from statistical physics literature \cite{jaynes2003probability} --- in addressing para-SDM problems; where in prior literature MEP has proven itself successful in addressing a wide variety of combinatorial optimization problems such as facility location problem \cite{rose1991deterministic}, protein structure alignment \cite{chen2005protein}, and graph and markov chain aggregation \cite{xu2014aggregation}. In brief, the MEP-based framework proposed in \cite{9517030} simultaneously determines a distribution over the {\em paths}, and the parameter values that maximize the associated Shannon entropy \cite{jaynes2003probability} of the distribution; while ensuring that the expected cumulative cost of the para-SDM attains a pre-specified value. The framework, then, employs an iterative scheme (an {\em annealing} scheme) and improves upon the distribution over the paths (i.e., the decision policy) as well as the parameter values that correspond to decreasing cumulative cost values of the SDM. In the current work we address the dynamic para-SDM problems where the parameters have an associated dynamics. These dynamic para-SDM problems belong to the class of combinatorial optimization problems where the underlying model parameters are time-varying. Prior literature such as \cite{sharma2012entropy, xu2013clustering, srivastava2020simultaneous} address specific instances of such time-varying combinatorial optimization problems. For instance, \cite{sharma2012entropy,xu2013clustering} address the data clustering problems where the underlying data points have associated velocity \cite{sharma2012entropy} or acceleration-driven \cite{xu2013distance} dynamics. Thus, requiring to determine the time varying clustering solutions as well as the dynamics of the resulting cluster centroids. The work done in \cite{srivastava2020simultaneous} addresses the dynamic facility location and route optimization problems in the context of spatial networks. However, the latter is restricted only to routes of pre-specified lengths (i.e., fixed finite horizon). Our current work is a generalization of the above problems to a much larger class of para-SDM problems, that feature optimization of paths over infinite horizons, and allow us to explicitly incorporate stochasticity and parameterization of the states, actions, cost function and dynamics underlying the para-SDMs.

One straightforward approach to address such problems is to solve the associated parameterized SDM at each time instant $t$, i.e. determine (a) the decision policy, and (b) the unknown parameters $\zeta(t)$ that minimizes the cumulative cost incurred by the para-SDM. However, there are several problems associated with such an approach. First, in many application areas the latter objective makes the para-SDM non-convex optimization problems as stated above. Solving these optimization problem repeatedly at every time instant is computationally expensive and unfit for application in an online environment. In fact, we observe computational times as large as $\mathcal{O}(10^4)$ times more than our proposed methodology (see Section \ref{sec: Simulations}); thereby making the latter a much more computationally viable approach under the assumption that the computational decision dynamics are faster in comparison to the SDM parameter dynamics. Second, under this frame-by-frame approach the resulting dynamics of the unknown parameters may not be viable; i.e., the {\em optimal values} of the unknown parameters $\Upsilon(t)$ and $\Upsilon(t+\Delta t)$ that minimize the cost function at time instants $t$ and $t+\Delta t$, respectively, could be {\em significantly} different from one another owing to the non-convexity, and thus, the dynamics that changes the parameter values from $\Upsilon(t)$ to $\Upsilon(t+\Delta t)$ becomes practically infeasible. We further elaborate on this in our simulations.

In this work we view the dynamic para-SDM problems from a control-theoretic viewpoint, and design control laws $u(\Upsilon,t)$ that aim at minimizing the instantaneous cumulative cost of the underlying para-SDM problem. The optimal decision policy of the SDM, which is also time-varying owing to the evolving parameters $\Upsilon(t)$, is evaluated as the fixed point of the recursive Bellman equation satisfied by underlying cumulative cost function. We build upon the Maximum Entropy Principle (MEP) based framework proposed in \cite{9517030} that addresses the static para-SDM problems. In particular, this MEP-based framework results into a {\em smooth approximation} (also referred to as {\em free-energy}) of the cumulative cost, which we exploit as a control-Lyapunov function describing the dynamic para-SDM problems. The main contributions here can be summarized as: (a) We formulate a non-linear feedback control law $u(\Upsilon,t)$ that governs the dynamics of the parameters $\Upsilon(t)$ in the para-SDM, and subsequently determine the time-varying optimal decision policy, (b) we show that under the proposed feedback control law, the parameter dynamics asymptotically tracks the local optimal of the underlying para-SDM problem (see Theorem \ref{thm: dynaMDPAsymptotic}), and (c) we show that this feedback control law is non-conservative, that is, if there
exists {\em a} Lipschitz control law that asymptotically tracks the local optimal of the para-SDM, then our proposed control law is also Lipschitz and bounded, and tracks a local optimal point (see Theorem \ref{thm: dynaMDPLipschitz}).

\section{Problem Formulation}\label{sec: PF}
We define a para-SDM as the tuple $\mathcal{M}=\langle \mathcal{S}, \mathcal{A},c,p,\gamma,\linebreak[1]\zeta,\linebreak[1]\eta\rangle$ where  $\mathcal{S}=\{s_1,\hdots,s_N=\delta\}$ denotes the state space with $s_N=\delta$ as a {\em cost-free termination} state, $\mathcal{A}=\{a_1,\hdots,a_M\}$ denotes the action space, $c:\mathcal{S}\times\mathcal{A}\times\mathcal{S}\rightarrow \mathbb{R}$ is the transition cost function; $p:\mathcal{S}\times\mathcal{S}\times\mathcal{A}\rightarrow [0,1]$ is the state transition probability function and $0 < \gamma \leq 1$ is a discounting factor; $\zeta=\{\zeta_s\in\mathbb{R}^{d_{\zeta}}:s\in\mathcal{S}\}$ and $\eta=\{\eta_a\in\mathbb{R}^{d_{\eta}}:s\in\mathcal{A}\}$ denote the state and action parameters, respectively. A decision policy $\mu:\mathcal{A}\times\mathcal{S}\rightarrow \{0,1\}$ determines the action taken at each state $s\in\mathcal{S}$, where $\mu(a|s) = 1$ implies that action $a\in\mathcal{A}$ is taken when the system is in the state $s\in\mathcal{S}$ and $\mu(a|s)=0$ indicates otherwise. For every initial state $x_0=s$ and (unknown) parameter values in $\zeta$ and $\eta$, the decision policy $\mu$ induces a stochastic process, whose realization is an infinite {\em path} $\omega$ (of consecutive actions and states) given by $\omega = (u_0, x_1, u_1, x_2, u_2,\hdots, x_K, u_K,x_{K+1},\hdots)$, where $u_k\in\mathcal{A}$, $x_k\in\mathcal{S}$ for all $k\in\mathbb{Z}_{\geq 0}$. The corresponding cumulative cost incurred is
\begin{align}\label{eq: ParValFunc}
J^{\mu}_{\zeta\eta}(s) = \mathbb{E}_{p_{\mu}}\Big[\sum_{k=0}^{\infty}\gamma^k c\big(x_k^{\zeta},u_k^{\eta},x_{k+1}^{\zeta}\big)\Big|x_0=s\Big],
\end{align}
where $x_k^{\zeta}$ denotes a state $x_k\in\mathcal{S}$ parameterized by $\zeta_{x_k}\in\zeta$, $u_k^{\eta}$ denotes the action $u_k\in\mathcal{A}$ parameterized by $\eta_{u_k}\in\eta$, and the expectation is with respect to the probability distribution $p_{\mu}(\cdot|s):\omega\rightarrow [0,1]$ on the space of all possible paths $\omega\in\Omega:=\{(u_k,x_{k+1})_{k\in\mathbb{Z}_{\geq 0}}: u_k\in\mathcal{A}, x_k\in\mathcal{S}\}$. To avoid notational clutter, we will drop the superscript in $x_k^{\zeta}$ and $u_k^{\eta}$ when clear from the context.

\begin{remark}
To ensure that the cumulative cost $J_{\zeta\eta}^{\mu}(s)$ is finite for all $s\in\mathcal{S}$ and the system reaches the cost-free termination state $\delta$ in finite steps, we assume that there exists atleast one {\em proper} policy $\bar{\mu}(a|s)\in\{0,1\}$ $\forall a\in\mathcal{A}$, $s\in\mathcal{S}$, and for all parameter values in $\zeta$ and $\eta$, under which there is a non-zero probability to reach the cost-free termination state $\delta$ starting from any state $s\in\mathcal{S}$ (please see \cite{9517030} for proof).
\end{remark}

In the case of {\em dynamic} para-SDM tasks, the parameter sets  $\zeta_1=\{\zeta_s\in\mathbb{R}^{d_\zeta}:s\in\mathcal{S}_1\subseteq\mathcal{S}\}$ and $\eta_1=\{\eta_a\in\mathbb{R}^{d_\eta}:a\in\mathcal{A}_1\subseteq\mathcal{A}\}$ denote the state and action parameters with {\em pre-specified} dynamics given by
\begin{align}\label{eq: givenDyna}
\dot{\zeta}_1 = \mathbf{\phi}_1(\zeta,\eta,t), \quad 
\dot{\eta}_1 = \mathbf{\psi}_1(\zeta,\eta,t),
\end{align}
where we assume that the dynamics $\phi_1,\psi_1\in C^1$ are continuously differentiable. Let $\zeta_2=\{\zeta_s:s\in\mathcal{S}_2=\mathcal{S}\text{\textbackslash}\mathcal{S}_1\}$ and $\eta_2=\{\eta_a:a\in\mathcal{A}_2=\mathcal{A}\text{\textbackslash}\mathcal{A}_1\}$ denote the parameters with {\em manipulable} dynamics. Owing to the time evolution of the parameters $\zeta_1$, $\eta_1$, the optimal decision policy $\mu^*$, and the state and action parameters $\zeta_2$, $\eta_2$ that minimize the cost function (\ref{eq: ParValFunc}) at each time instant $t$ are also time-varying. More precisely, the objective is to determine the evolution of the parameters $\zeta_2$, $\eta_2$ and the decision policy $\mu$ such that the cumulative cost 
\begin{align}\label{eq: dynamic_ParValFunc}
J^{\mu}_{\Upsilon}(s,t) = \mathbb{E}_{p_{\mu}}\Big[\sum_{k=0}^{\infty}\gamma^k c\big(x_k^{\zeta(t)},u_k^{\eta(t)},x_{k+1}^{\zeta(t)}\big)\big|x_0=s\Big]
\end{align}
is minimized at for all $t\in\mathbb{R}_{\geq 0}$, where $\Upsilon=[\zeta_1,\eta_1,\zeta_2,\eta_2]$. We propose a control-theoretic framework to determine
\begin{align}\label{eq: findDyna}
\dot{\zeta}_2 = \mathbf{\phi}_2(\zeta,\eta,t),\quad 
\dot{\eta}_2 = \mathbf{\psi}_2(\zeta,\eta,t)
\end{align}
dynamics that governs the time evolution of the parameters $\zeta_2$, and $\eta_2$, and subsequently evaluate the decision policy $\mu^*$ at each time instant $t$. As the work done in \cite{9517030}, that addresses static para-SDM, forms the foundation for our control-theoretic framework to address dynamic para-SDM, we briefly illustrate it in the next section.

\section{MEP-based Approach to {\em Static} Para-SDM}\label{sec: MEP_paraSDM}
MEP states that the {\em most unbiased} probability distribution $p_{\mathcal{X}}(\cdot)$ of a random variable $\mathcal{X}$ under the constraint on the expected value of the functions $\nu_k:\mathcal{X}\rightarrow \mathbb{R}$ for all $1\leq k\leq m$ is the one that solves
\begin{align}
\begin{split}
\max_{\{p_{\mathcal{X}}(x_i)\}}\quad &\mathcal{H}(\mathcal{X}) = -\sum_{i=1}^n p_{\mathcal{X}}(x_i)\ln p_{\mathcal{X}}(x_i)\\
\text{subject to}\quad &\sum_{i=1}^n p_{\mathcal{X}}(x_i)\nu_k(x_i) = N_k~\forall~1\leq k\leq m,
\end{split}
\end{align}
where $N_k$ for all $1\leq k\leq m$ are given. The framework proposed in \cite{9517030} employs MEP to address the para-SDM problems. More precisely, the framework determines the {\em most unbiased} distribution $\{p_{\mu}(\omega|s):\omega\in\Omega\}$ such that the cost function $J_{\Upsilon}^{\mu}(s)$ attains a given value $J_0$, i.e. it solves the following optimization problem
\begin{align}\label{eq: OptimP1}
\begin{split}
\max_{\{p_{\mu}(\cdot|s)\}:\mu\in\Gamma}\quad \mathcal{H}^{\mu}_{\Upsilon}(s) &= -\sum_{\omega\in\Omega}p_{\mu}(\omega|s)\log p_{\mu}(\omega|s)\\
\text{subject to}\quad J^{\mu}_{\Upsilon}(s) &= J_0,
\end{split}
\end{align}
where $\Gamma:=\{\pi:0<\rho\leq \pi(a|s)<1\forall a\in\mathcal{A},s\in\mathcal{S}\}$ denotes the set of all {\em proper} control policies, and $\rho>0$ is arbitrarily small. The optimization problem (\ref{eq: OptimP1}) is well-posed since the maximum entropy $\mathcal{H}^{\mu}_{\Upsilon}(s)$ for all $s\in\mathcal{S}$ is finite for the class of proper policy $\mu\in\Gamma$ \cite{biondi2014maximizing,savas2018entropy}. Note that the policy $\mu$ is defined over finite action and state spaces, whereas the $p_{\mu}(\omega|s)$ is defined over infinitely many paths $\omega\in\Omega$; thus, \cite{9517030} exploits the Markov property that dissociates $p_{\mu}(\omega|s)$ in terms of the policy $\mu$ and state transition probability as $p_{\mu}(\omega|x_0) = \prod_{t=0}^{\infty}\mu(u_k|x_k)p(x_{t+1}|x_k,u_k),$ and prudently chooses to work with $\mu$ instead of $p_{\mu}$. The Lagrangian for (\ref{eq: OptimP1}) is given by $V_{\beta\Upsilon}^{\mu}(s):=J^{\mu}_{\Upsilon}(s)-\frac{1}{\beta}\mathcal{H}^{\mu}(s)$ \footnote{The Lagrange parameter $\beta$ decides the constraint value $J_0$ in (\ref{eq: OptimP1}). Thus, for a given $\beta$, $J_0$ is a constant and we ignore it in the expression of the Lagrangian $V_{\beta\Upsilon}^{\mu}(s)$.}, that follows the following recursive Bellman equation
\begin{align}\label{eq: BellmanTrue}
V_{\beta\Upsilon}^{\mu}(s) = \sum_{\substack{a,s'}}\mu_{a|s}p_{ss'}^{a}\big(\bar{c}_{ss'}^{a,\mu}+\gamma V_{\beta\Upsilon}^{\mu}(s')\big)+c_0(s),
\end{align}
where $\mu_{a|s}=\mu(a|s)$, $p_{ss'}^a=p(s'|s,a)$,  $\bar{c}_{ss'}^{a,\mu}=c(s,a,s')+\frac{\gamma}{\beta}\log p(s'|s,a)+\frac{\gamma}{\beta}\log\mu_{a|s}$ for simplicity in notation, and $c_0(s)$ depends on $\gamma$ and $\beta$, and  is independent of the policy $\mu$ and the parameters $\Upsilon$. Without loss of generality, we ignore $c_0(s)$ in the subsequent calculations (see \cite{srivastava2021paraSDM}). For proof of the above Bellman equation please see Theorem 1 in \cite{9517030} (or detailed proof in \cite{srivastava2021paraSDM}). The optimal policy $\mu_{\beta}^*$ is obtained by setting $\frac{\partial V_{\beta\Upsilon}^{\mu}(s)}{\partial \mu(a|s)}=0$, which results into the Gibbs' distribution
\begin{align}
&\mu^*_{\beta\Upsilon}(a|s) = \frac{\exp\big\{-(\beta/\gamma) \Lambda_{\beta\Upsilon}(s,a)\big\}}{\sum_{a'\in\mathcal{A}}\exp\big\{-(\beta/\gamma)\Lambda_{\beta\Upsilon}(s,a')\big\}},\label{eq: Policy}\\
&\text{where }\Lambda_{\beta\Upsilon}(s,a) = \sum_{s'\in\mathcal{S}}p_{ss'}^{a}\big(\bar{c}_{ss'}^{a} + \gamma V^{*}_{\beta\Upsilon}(s')\big)\label{eq: Q_bell}
\end{align}
is the state-action value function, 
$p_{ss'}^a = p(s'|s,a)$, $c_{ss'}^a = c(s,a,s')$,  $\bar{c}_{ss'}^a=c_{ss'}^a+\frac{\gamma}{\beta}\log p_{ss'}^a$, and
\begin{align}\label{eq: free_energy}
V^{*}_{\beta\Upsilon}(s) = -\frac{\gamma}{\beta}\log\Big(\sum_{a\in\mathcal{A}}\exp\Big\{-\frac{\beta}{\gamma}\Lambda_{\beta\Upsilon}(s,a)\Big\}\Big)
\end{align}
is the value function corresponding to the global optimal policy $\mu^*_{\beta\Upsilon}$ that is obtained on substituting (\ref{eq: Policy}) in (\ref{eq: BellmanTrue}). Note that from (\ref{eq: Q_bell}) and (\ref{eq: free_energy}) it can be deduced that the state-action value function $\Lambda_{\beta\Upsilon}(s,a)$ satisfies the implicit equation $\Lambda_{\beta\Upsilon}(s,a)=:[T\Lambda_{\beta\Upsilon}](s,a)$, where
\begin{align}\label{eq: Q_map}
&[T\Lambda_{\beta\Upsilon}](s,a) = \sum_{s'\in\mathcal{S}}p_{ss'}^{a}\big(c_{ss'}^{a} + \frac{\gamma}{\beta}\log p_{ss'}^{a}\big)\nonumber\\
&\quad -\frac{\gamma^2}{\beta}\sum_{s'\in\mathcal{S}}p_{ss'}^{a}\log\sum_{a'\in\mathcal{A}}\exp\Big\{-\frac{\beta}{\gamma}\Lambda_{\beta\Upsilon}(s',a')\Big\}
\end{align}
is a contraction map (see Theorem 2 in \cite{9517030}), and thus, $\Lambda_{\beta\Upsilon}$ is a fixed-point of the map $T$. Subsequently, the unknown state and action parameters $\zeta_2$ and $\eta_2$ are determined by setting $\sum_{s'\in\mathcal{S}}\frac{\partial V_{\beta\Upsilon}^*(s')}{\partial \zeta_s}=0\text{ and }\sum_{s'\in\mathcal{S}}\frac{\partial V_{\beta\Upsilon}^*(s')}{\partial \eta_a}=0$ for all $s\in\mathcal{S}_2$ and $a\in\mathcal{A}_2$. It is straightforward from (\ref{eq: free_energy}) that the partial derivative $\frac{\partial V_{\beta\Upsilon}^*(s')}{\partial \zeta_s}=:G_{\zeta_s}^\beta(s')$ satisfy the recursive Bellman
\begin{align}\label{eq: derivBell}
G_{\zeta_s}^\beta(s') = \sum_{a',s''}\mu_{\beta}(a'|s')p_{s's''}^{a'}\Big[\frac{\partial c_{s's''}^{a'}}{\partial \zeta_s}+\gamma G_{\zeta_s}^\beta(s'')\Big]
\end{align}
The recursive Bellman equation for $\frac{\partial V_{\beta\Upsilon}^*(s')}{\partial \eta_a}=:G_{\eta_a}^{\beta}(s')$ is similarly derived.

\section{A Control Theoretic Framework for Dynamic Parameterized SDM}\label{sec: PS}
As briefly stated in the Section \ref{sec: Intro}, a straightforward method to determine the dynamical evolution of the parameters $\zeta_2(t)$ and $\eta_2(t)$ would be to solve the optimization problem (\ref{eq: ParValFunc}) at each time instant $t$. However, such a method is computationally expensive, and does not make use of the {\em known} dynamical evolution of the parameters $\zeta_1$ and $\eta_1$; resulting into much computational redundancy. Additionally, minimizing the (possible non-convex) cost function $J^\mu_{\Upsilon}(s,t)$ in (\ref{eq: dynamic_ParValFunc}) may result into (local) optimal values of $\zeta_2(t),\eta_2(t)$ and $\zeta_2(t+\Delta t),\eta_2(t+\Delta t)$ of the unknown parameters that are {\em significantly} different from each other; thereby, resulting into a non-viable dynamics (\ref{eq: findDyna}). In this work, we build up on the MEP-based solution to the static para-SDM problem (in Section \ref{sec: MEP_paraSDM}), and propose a control-theoretic framework that addresses the above issues. In particular, we exploit the use of free-energy $V_{\beta}^*(s)$ in (\ref{eq: free_energy}) in determining an appropriate control-Lyapunov candidate for the dynamical systems (\ref{eq: givenDyna}) and (\ref{eq: findDyna}) represented together as the control system
\begin{align}\label{eq: CtrlSystm}
\dot{\Upsilon}=f(\Upsilon,u,t),
\end{align}
where $\Upsilon=[\zeta_1~\eta_1~\zeta_2~\eta_2]\in\mathbb{R}^{Nd_{\zeta}+Md_{\eta}}$, $f=[\phi_1~\psi_1~\phi_2~\psi_2]\in\mathbb{R}^{Nd_{\zeta}+Md_{\eta}}$. Our objective is to design the control field $u(t):=[\phi_2~\psi_2]$ such that the value function (\ref{eq: dynamic_ParValFunc}) is minimized at each time instant. We consider the following control-Lyapunov candidate
\begin{align}\label{eq: ctrl_Lyap}
\text{$\mathcal{V}(\Upsilon)=\sum_{s\in\mathcal{S}}V_{\beta\Upsilon}^*(s),$}
\end{align}
where $V_{\beta\Upsilon}^*(s)$ is the free energy given in (\ref{eq: free_energy}) that corresponds to the optimal decision policy $\mu^*_{\beta\Upsilon}$ in (\ref{eq: Policy}). Subsequently, we determine the control field $u(t)$ (i.e., the dynamics $\phi_2$ and $\psi_2$) such that the time derivative $\dot{\mathcal{V}}$ along the trajectory $\Upsilon(t)$ is non-positive.\\
\begin{remark}
The above choice of Lyanpunov candidate $V(\Upsilon)$ is essential in being able to design a control law $u(t)$ that governs the evolution of the parameters $\zeta_2$ and $\eta_2$. In particular, $V_{\beta\Upsilon}^*(s)$ in (\ref{eq: ctrl_Lyap}) is a smooth approximation of the cumulative costs $J_{\Upsilon}^{\mu^*}(s)=\min_{\mu}J_{\Upsilon}^{\mu}(s)$ at the optimal decision policy $\mu^*:\mathcal{S}\times\mathcal{A}\rightarrow\{0,1\}$, and is only a function of the state and action parameters $\Upsilon$. As illustrated shortly, the time derivative $\dot{\mathcal{V}}(\Upsilon):=\frac{\partial \mathcal{V}}{\partial \Upsilon}\frac{d\Upsilon}{dt}$ is an affine function of the control law $u(t)$ --- that makes it easy to design appropriate $u(t)$ such that $\dot{\mathcal{V}}(\Upsilon)$ is non-negative.
\end{remark}

We further summarize the properties of $\mathcal{V}$, and its time derivative in the following theorem.

\begin{theorem}\label{thm: DynaMDP_Lyap}
Let $\mathcal{V}$ be the control-Lyapunov candidate for the dynamical system given by (\ref{eq: CtrlSystm}). Then,\\
(a) Positive definiteness: There exists a constant $c>0$ such that $\mathcal{V}(\Upsilon) + c>0$ for all $\Upsilon\in\mathbb{R}^{Nd_{\zeta}+Md_{\eta}}$.\\
(b) There is no dynamic control authority only at time instants when the unknown parameters $\zeta_2,\eta_2$ are at local minima.
\end{theorem}
{\em Proof. }(a) We have that $\mathcal{V}(\Upsilon)=\sum_{s\in\mathcal{S}}V_{\beta\Upsilon}^*(s)=\sum_{s\in\mathcal{S}}\min_{\mu\in\Gamma}V_{\beta\Upsilon}^{\mu}(s)$, where $V_{\beta\Upsilon}^{\mu}(s)$ is the Lagrangian for the optimization problem in (\ref{eq: OptimP1}). Thus, $\mathcal{V}(\Upsilon)\geq \sum_{s\in\mathcal{S}}\big(\min_{\mu\in\Gamma}J^\mu_{\Upsilon}(s)-\frac{1}{\beta}\max_{\mu\in\Gamma}\mathcal{H}^{\mu}(s)\big)$. The entropy $\mathcal{H}^{\mu}(s)<\infty$ for all $\mu\in\Gamma$, and $\min_{\mu\in\Gamma}J^{\mu}_{\Upsilon}(s)>0$. Thus, there exists a $c>0$ such that $\mathcal{V}\geq -c$.\\
(b)The time-derivative $\dot{\mathcal{V}}$ is given by
\begin{align}\label{eq: timeDerivDynaMDP}
\dot{\mathcal{V}}=G^T\kappa+F^Tu,
\end{align}
where $\kappa^T=[\phi_1^T~\psi_1^T]$, $G^T=[G_\phi^T~G_\psi^T]\in\mathbb{R}^{|\mathcal{S}_1|d_{\zeta}+|\mathcal{A}_1|d_{\eta}}$, $G_\phi=[g_\phi(s)]_{s\in\mathcal{S}_1}\in\mathbb{R}^{|\mathcal{S}_1|d_{\zeta}}$, $g_\phi(s)=\frac{\partial \mathcal{V}(\Upsilon)}{\partial \zeta_s}$, $G_\psi=[g_\psi(a)]_{a\in\mathcal{A}_1}\in\mathbb{R}^{|\mathcal{A}_1|d_{\eta}}$, $g_\psi(a)=\frac{\partial \mathcal{V}(\Upsilon)}{\partial \eta_a}$, $F^T=[F_\phi^T~F_\psi^T]\in\mathbb{R}^{|\mathcal{S}_2|d_{\zeta}+|\mathcal{A}_2|d_\eta}$, $F_\phi=[f_\phi(s)]_{s\in\mathcal{S}_2}\in\mathbb{R}^{|\mathcal{S}_2|d_{\zeta}}$, $f_\phi(s)=\frac{\partial \mathcal{V}(\Upsilon)}{\partial \zeta_s}$, $F_\psi=[f_\psi(a)]_{a\in\mathcal{A}_2}\in\mathbb{R}^{|\mathcal{A}_2|d_{\eta}}$, $f_\psi(a)=\frac{\partial \mathcal{V}(\Upsilon)}{\partial \eta_a}$. There is no dynamic control authority when $\frac{\partial \dot{\mathcal{V}}}{\partial u}=F$ is zero, i.e., $F=0$. By above definition of $F$, we have that when F=0, $\frac{\partial \mathcal{V}}{\partial \zeta_s}=0$ $\forall$ $s\in\mathcal{S}_2$ and $\frac{\partial \mathcal{V}}{\partial \eta_a}=0$ $\forall$ $a\in\mathcal{A}_2$.

{\em Control Design for tracking parameters $\zeta_2,\eta_2$: } We make use of the {\em affine} dependence of $\dot{\mathcal{V}}$ in (\ref{eq: timeDerivDynaMDP}) on the control $u(t)$ to make $\dot{\mathcal{V}}$ non-positive analogous to the chosen control design \cite{sontag1983lyapunov,sepulchre2012constructive}. Specifically, we choose the control of the form
\begin{align}\label{eq:  Control_choice1}
\text{\small$u(F) = -\Bigg[K_0 + \dfrac{\alpha + \sqrt{\alpha^2+(F^TF)^2}}{F^TF}\Bigg]F,$}
\end{align}
when $F\neq 0$, and $u(F)=0$ otherwise; here $\alpha = G^T\kappa$, and $K_0>0$. The following theorems show that given the dynamics (\ref{eq: givenDyna}) for the parameters $\zeta_1$, $\eta_1$ to be continuously differentiable, the state and action parameters $\zeta_2(t),\eta_2(t)$ asymptotically track the condition $\|F\|_2\rightarrow 0$, where $\|\cdot\|_2$ indicates $2$-norm. Further, if there exists a bounded control design $\hat{u}(t)$ that ensures $\dot{V}\leq 0$, then our proposed control $u(F)$ in (\ref{eq:  Control_choice1}) is also bounded, i.e., our control design $u(F)$ is not conservative.

\begin{theorem}\label{thm: dynaMDPAsymptotic}
Asymptotic convergence: For the dynamical system (\ref{eq: CtrlSystm}) the choice of control $u(F)$ in (\ref{eq:  Control_choice1}) results in $\dot{\mathcal{V}}\leq 0$ $\forall$ $t \geq 0$ and the derivatives in $\|F(t)\|_2 \rightarrow 0$ as $t\rightarrow \infty$; i.e., the state and action parameters $\zeta_2(t)$ and $\eta_2(t)$ asymptotically track the local optimal.
\end{theorem}
{\em Proof.} Substituting the control law $u(F)$ in (\ref{eq: Control_choice1}) into the expression of $\dot{\mathcal{V}}$ in (\ref{eq: timeDerivDynaMDP}) we obtain $\dot{\mathcal{V}} = -K_0F^TF-\sqrt{(G^T\kappa)^2+(F^TF)^2}$, where $K_0>0$. Hence $\dot{\mathcal{V}}\leq 0$. We also know from Theorem \ref{thm: DynaMDP_Lyap} that the function $\mathcal{V}$ is lower bounded. Thus, $\mathcal{V}(t)$ converges to $\mathcal{V}_{\infty}$, where $|\mathcal{V}_{\infty}|<\infty$.

Note that $\int_{0}^{\infty}|\dot{\mathcal{V}}(\tau)|d\tau=\mathcal{V}_{\infty}-\mathcal{V}(0)\leq 0$, and under the assumption that the dynamics (\ref{eq: givenDyna}) of the parameters $\zeta_1$ and $\eta_1$ are continuous and differentiable, $\dot{\mathcal{V}}$ is of bounded variation \cite{royden1988real}. Thus, by Lemma \ref{lem: boundedVar} (see Appendix) we have $|\dot{\mathcal{V}}|\rightarrow 0$ as $t\rightarrow \infty$. Now since $\dot{\mathcal{V}} = -K_0F^TF-\sqrt{(G^T\kappa)^2+(F^TF)^2}$, we have that $K_0F^TF \leq |\dot{\mathcal{V}}|$. Thus we conclude that $\|F(t)\|_2 \rightarrow 0$ as $t\rightarrow \infty$.

\begin{theorem}\label{thm: dynaMDPLipschitz}
Lipschitz continuity: If there exists a control $\hat{u}:\mathbb{R}^{Nd_{\zeta}+Md_{\eta}}\rightarrow \mathbb{R}^{|\mathcal{S}_2|d_{\zeta}+|\mathcal{A}_2|d_{\eta}}$ Lipschitz at $F=0$, such that $\dot{\mathcal{V}}\leq 0$ $\forall$ $t\geq 0$ for this control, then the choice of control $u(F)$ in (\ref{eq:  Control_choice1}) is Lipschitz at $F=0$. That is, $\exists$ $\epsilon > 0$ and a $c_0$ such that $\|u(F)\|\leq c_0\|F\|$ for $\|F\|\leq \epsilon$.
\end{theorem}
{\em Proof.} Note - the proof here is similar to the proof for the Proposition 3.43 in \cite{sepulchre2012constructive}. Since $\hat{u}$ is Lipschitz at $F=0$, there exists a neighbourhood $B_r=\{F:\|F\|\leq r\}$ and $\bar{k}>0$ such that $\|\hat{u}\|\leq \bar{k}\|F\|$ for all $F\in B_r$. By definition, $\dot{\mathcal{V}}=G^T\kappa+F^T\hat{u}\leq 0$. If $G^T\kappa>0$, then $|G^T\kappa|\leq|F^T\hat{u}|\leq\|F\|\|\hat{u}\|\leq\bar{k}\|F\|\|F\|$ $\forall$ $F\in B_r$. Thus, the control design $u(F)$ in (\ref{eq:  Control_choice1}) can be bounded above as $\|u(F)\|\leq (K_0+2\bar{k}+1)\|F\|$. For the case when $G^T\kappa<0$, we have that $\|u(F)\|\leq (1+K_0)\|F\|\leq (1+K_0)\|F\|$.

\begin{figure*}
\centering
\includegraphics[width=0.95\textwidth]{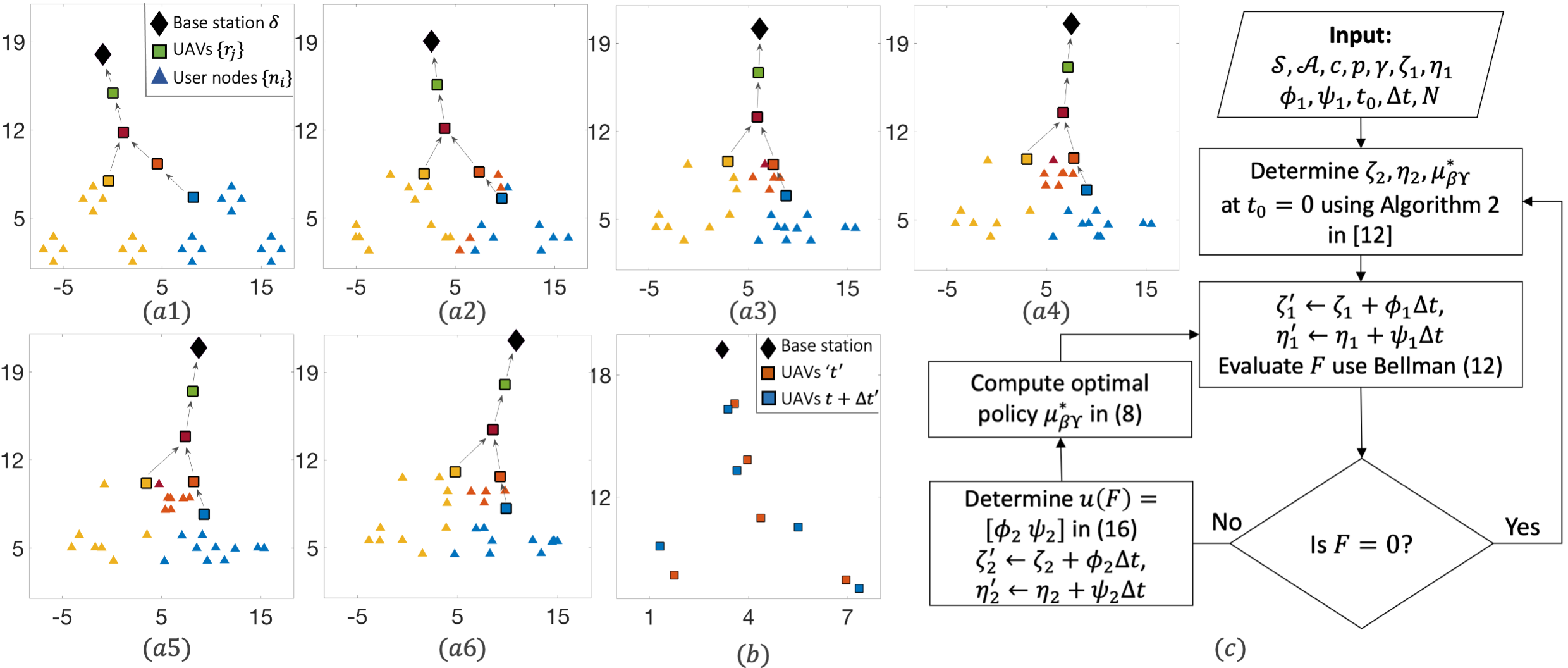}
\caption{(a1)-(a6) Illustrates the dynamic multi-UAV network problem. Observe the change in locations of the user nodes $\{n_i\}$ and the aerial base station $\delta$ with time. Thus, the resulting UAV locations in the network have dynamics governed by $u(t)$ in (\ref{eq:  Control_choice1}). Also, observe the change in the color of the triangles (denoting user nodes) from (a1) to (a2), (a2) to (a3) (, and further) indicating the change in communication paths. (b) Indicates an drastic change in UAV locations at time $t$ and $t+\Delta t$. (c) Algorithm for dynamic para-SDM.}
\label{fig: dynaMDP_simulation}
\end{figure*}

{\em Algorithmic Insights: } As stated in the Section \ref{sec: PF}, we need to determine the dynamical evolution of the parameters $\zeta_2$ and $\eta_2$, as well as the time-varying optimal policy $\mu_{\beta}^*$ such that the cumulative cost $J_{\Upsilon}^{\mu}(s,t)$ is minimized at each time instant $t$ for all $s\in\mathcal{S}$. In our proposed methodology the control law $u(F)$ in (\ref{eq:  Control_choice1}) addresses the dynamical evolution of the parameters $\zeta_2$ and $\eta_2$ only. We do not explicitly design a control law governing the time-varying optimal policy $\mu_{\beta\Upsilon}^*$ of the para-MDP, albeit we directly make use of the expression in (\ref{eq: Policy}) to determine the optimal policy at each time instant. This requires solving for the fixed point $\Lambda_{\beta\Upsilon}(s,a)$ of the contraction map in (\ref{eq: Q_map}). In case of limited computational resources (or when solving for the fixed point is computationally intensive), we alternatively propose to use the expression in (\ref{eq: Q_bell}) to estimate $\Lambda_{\beta\Upsilon}(s,a)$. In particular, we consider the first order Taylor series approximation of value function $V_{\beta\Upsilon}^*(s)$ given by
\begin{align}\label{eq: taylorApprox}
\text{\small $V_{\beta\Upsilon'}^*(s) \approx V_{\beta\Upsilon}^*(s) + \sum_{s',a'}\frac{\partial V_{\beta\Upsilon}^*(s)}{\partial \zeta_{s'}}\delta\zeta_{s'}+\frac{\partial V_{\beta\Upsilon}^*(s)}{\partial \eta_{a'}}\delta\eta_{a'}$}
\end{align} 
where $\Upsilon'=\Upsilon+\delta\Upsilon$, to approximate the value function $V_{\beta\Upsilon'}^*(s)$ at time $t+\Delta t$ (i.e., when the parameters are $\Upsilon'$) using the known optimal value function $V_{\beta\Upsilon}^*(s)$ at time $t$ (i.e., when the parameters are $\Upsilon$). Subsequently, the fixed point $\Lambda_{\beta\Upsilon}(s,a)$ is estimated as in (\ref{eq: Q_bell}), and the optimal policy $\mu_{\beta}^*$ is given by (\ref{eq: Policy}). Please refer to the Algorithm in Figure \ref{fig: dynaMDP_simulation}(c) for detailed steps.
\section{Simulations}\label{sec: Simulations}
In this section we simulate a dynamic para-SDM problem to demonstrate the efficacy of our proposed control design $u(F)$ in (\ref{eq:  Control_choice1}). In particular, we consider the multi-UAV network systems (see Figure \ref{fig: 5Gsmall_cell}(b)) illustrated in the Section \ref{sec: Intro}. In such networks, the user nodes and base station are usually mobile, i.e., the locations $\{x_i\in\mathbb{R}^d\}$ and $z\in\mathbb{R}^d$ of the user nodes $\{n_i\}$ and the aerial base station $\delta$, respectively, are time-varying. The objective is to determine the dynamical evolution (of the locations $\{y_j\in\mathbb{R}^d\}$) of the UAVs $\{r_j\}$, and the time-varying multi-hop communication path (via the network of UAVs) from each user $n_i$ to the base station $\delta$ such that the total communication cost gets minimized at each time instant.

We model the multi-UAV network as a para-SDM $\mathcal{M}=\langle \mathcal{S},\mathcal{A},c,p,\zeta\rangle$. Here, the state space $\mathcal{S}=\{\{n_i\},\{r_j\},\delta\}$ comprises of all the user nodes, UAVs, and the base station, the action space $\mathcal{A}=\{\{r_j\},\delta\}$ consists of the UAVs and the base station such that any action $a\in\mathcal{A}$ indicates a communication packet hop at $a$, the state parameters $\zeta := \zeta_1\sqcup \zeta_2$ where $\zeta_1=\{\{x_i\},z\}$ is the set of parameters (user nodes and base station locations) with {\em predefined} (fixed) dynamics, and $\zeta_2=\{\{y_j\}\}$ is the set of parameters (the UAV locations) with {\em manipulable} dynamics, the cost function $c(s,a,s')=\|\zeta_s-\zeta_{s'}\|^2$ is the squared-euclidean distance between the spatial locations $\zeta_s\in\mathbb{R}^d$ and $\zeta_{s'}\in\mathbb{R}^d$ of the states $s$ and $s'$, respectively, and the transition probability $p(s'|s,a) = 1$ if $s'=a$ and $p(s'|s,a)=0$ otherwise for all $s,s'\in\mathcal{S}$ and $a\in\mathcal{A}$. Figure \ref{fig: dynaMDP_simulation}(a) illustrates the user nodes $\{n_i\}$, the aerial base station $\delta$ along with the UAVs $\{r_j\}$ (located at $\{y_j\}$), and the optimal communication paths at the time instant $t=0$. As illustrated in the Figure \ref{fig: dynaMDP_simulation}(a1), a user node $n_i$ of a particular color first sends its information packet to the UAV of the similar color which then reaches the base station via the indicated path. Note that the UAV locations and the communication paths at $t=0$ are as obtained using the Algorithm 2 in \cite{9517030} that addresses the static para-SDM scenario (i.e. when the user nodes and base station are considered stationary).

For the purpose of simulation, we assign randomly generated dynamics to the locations of the user nodes $\{n_i\}$ and the base station $\delta$ (i.e., to all the state parameters in $\zeta_1$) in the multi-UAV network; their corresponding spatial evolution is noted in the Figures \ref{fig: dynaMDP_simulation}(a2)-\ref{fig: dynaMDP_simulation}(a6). Please refer to the supplementary video material for a detailed illustration. We use the control design $u(F)$ proposed in (\ref{eq:  Control_choice1}) to determine the dynamical evolution of the state parameters in $\zeta_2$ (i.e., the time-varying UAV locations $\{y_j\}$); the corresponding communication paths are governed by the optimal policy $\mu_{\beta\Upsilon}^*$ in (\ref{eq: Policy}), i.e. $\mu_{\beta\Upsilon}^*(a|s)=1$ indicates that the communication packet at $s\in\mathcal{S}$ goes next to the UAV indicated by $a\in\mathcal{A}$. As illustrated in the Figure \ref{fig: dynaMDP_simulation}(a1), at time $t=0$ all the user nodes are coloured either blue or dark yellow, where the blue (dark yellow) user nodes send their communication packet to the blue (dark yellow) UAV; the subsequent communication path is as indicated by the arrows in the figure. As time progresses (see Figures \ref{fig: dynaMDP_simulation}(a2)-\ref{fig: dynaMDP_simulation}(a6)), the locations of the user nodes and the base station evolve based on their respective inherent dynamics, and the UAV locations $\{y_j\}$ evolve as per the control law $u(F)$ in (\ref{eq:  Control_choice1}). The corresponding time-varying communication paths from each user node $n_i$ to the base station $\delta$ is clearly indicated by the change in the color of the triangles representing the user nodes. Please refer to the supplementary video material for more details.

As briefly illustrated in the Section \ref{sec: Intro}, a straightforward way to determine the dynamical evolution of the UAV locations $\{y_j\}$ and the time-varying communication paths is to solve the optimization problem in (\ref{eq: dynamic_ParValFunc}) at each time instant $t$. However such an approach has apparent downsides to it. For instance, it is computationally intensive and requires approximately $\mathcal{O}(10^4)$ times more computational time to determine the evolution of the multi-UAV network illustrated in the Figure \ref{fig: dynaMDP_simulation}. Secondly, the underlying optimization problem is non-convex in nature (as it involves allocating UAVs in the network which is analogous to the non-convex Facility location problem \cite{mahajan2009planar}). Thus, the UAV locations $\{y_j\}$ obtained at time $t$ and $t+\Delta t$ could possibly be far from one another; resulting into a {\em non-viable} dynamics of the UAVs. For instance, the Figure \ref{fig: dynaMDP_simulation}(b) illustrates the UAV locations obtained at times $t$ and $t+\Delta t$ (where $\Delta t=0.1$ units). Note that there is a considerable change in the UAV locations at the two time instants that are only $\Delta t$ units apart; the control effort required to implement such non-viable dynamics will be vast, and possibly impractical.

\section{Conclusion}
This work develops a control-theoretic framework to address the class of time-varying para-SDM problems. These optimization problems require simultaneously determining the manipulable parameter dynamics as well as the time-varying policy such that the associated cumulative cost is minimized at each time instant. The ontrol design methodology presented in this work generalizes to the parameterized SDM problems with additional capacity or exclusion constraints on the state and action spaces. For instance, in the multi-UAV network systems in Figure \ref{fig: 5Gsmall_cell}(b) each UAV $r_j$ may be of limited capacity $c_j\in\mathbb{R}$ in terms of the number of users $\{n_i\}$ that it covers, i.e., $\sum_{s,a}p(r_j|s,a)\mu^*_{\beta\Upsilon}(a|s)\leq c_j$ \cite{9517030}. Such constraints reflect only in the optimal decision policy $\mu_{\beta\Upsilon}^*(a|s)$, and thus can easily be incorporated in our proposed scheme. The current work also extends to the class of problems where the transition costs $c(x_k^{\zeta},u_k^{\eta},x_{k+1}^{\zeta},\Gamma)$ as well as the state transition probability $p(x_{k+1}^{\zeta}|x_k^{\zeta},u_k^{\eta},\Xi)$, of the underlying SDM, are explicitly parameterized by $\Gamma$ and $\Xi$, respectively; for instance, these parameters could account for the effect of ambience and weather on the costs and signal transmission in case multi-UAV networks. 

\appendix
\begin{lemma}\label{lem: boundedVar}
Given that $\theta:\mathbb{R}\rightarrow\mathbb{R}_{\geq 0}$ is of bounded variation such that $\int_0^{\infty}\theta(\tau)d\tau<\infty$, then $\lim_{t\rightarrow\infty}\theta(t)=0$.
\end{lemma}
{\em Proof: }
By contradiction. Let $\lim_{t\rightarrow\infty}\theta(t)\neq 0$, then there exists an $\epsilon>0$ and a sequence $(x_n)$ of real numbers such that $\theta(x_n)>\epsilon$ and $x_{n+1}-x_n > 1$ for all $n\in\mathbb{N}$. Consider the function $g(t):=\int_0^{t}\theta(\tau)d\tau$ --- a monotonically increasing function bounded above by $\lambda:=\int_0^{\infty}\theta(\tau)d\tau$ --- that converges to $\lambda$ as $t\rightarrow\infty$. By definition, for every $\epsilon'>0$ there exists a $X_{\epsilon'}\in\mathbb{R}$ such that $\lambda - g(X_{\epsilon'})<\epsilon'$. More precisely, there exists a subsequence $(x_{n_k})$ such that $\int_{x_{n_k}}^{\infty}\theta(\tau)d\tau<(0.5)^{n_k}\epsilon$. Let $\theta(r_{n_k}):=\min_{\tau\in[x_{n_k},x_{n_{k+1}}]}\theta(\tau)$. Thus, $\theta(r_{n_k})<\theta(r_{n_k})(x_{n_{k+1}}-x_{n_k})\leq \int_{x_{n_k}}^{x_{n_{k+1}}}\theta(\tau)d\tau\leq\int_{x_{n_k}}^{\infty}\theta(\tau)d\tau<(0.5)^{n_k}\epsilon$. The variation of $\theta(t)$ on the interval $\mathcal{I}:=\sqcup_{k=1}^m(x_{n_k},r_{n_k})$ is given by $V_{\mathcal{I}}(\theta) = \sum_{k=1}^{m}|\theta(r_{n_k})-\theta(x_{n_k})|\geq \sum_{k=1}^{m}\big(\theta(x_{n_k}) - \theta(r_{n_k})\big)\geq m\epsilon - \sum_{n_k=0}^{\infty}(0.5)^{n_k}\epsilon = (m-1)\epsilon$. Therefore, the total variation of $\theta(t)$ is not finite which is a contradiction. Thus, $\lim_{t\rightarrow\infty}\theta(t)=0$.

\bibliographystyle{unsrtnat}        
\bibliography{autosam}           

\begin{thebibliography}{26}
\providecommand{\natexlab}[1]{#1}
\providecommand{\url}[1]{\texttt{#1}}
\expandafter\ifx\csname urlstyle\endcsname\relax
  \providecommand{\doi}[1]{doi: #1}\else
  \providecommand{\doi}{doi: \begingroup \urlstyle{rm}\Url}\fi

\bibitem[Littman(1996)]{littman1996algorithms}
Michael~Lederman Littman.
\newblock \emph{Algorithms for sequential decision-making}.
\newblock Brown University, 1996.

\bibitem[Liberzon(2011)]{liberzon2011calculus}
Daniel Liberzon.
\newblock \emph{Calculus of variations and optimal control theory}.
\newblock Princeton university press, 2011.

\bibitem[Bertsekas and Tsitsiklis(1996)]{bertsekas1996neuro}
Dimitri~P Bertsekas and John~N Tsitsiklis.
\newblock \emph{Neuro-dynamic programming}, volume~5.
\newblock Athena Scientific Belmont, MA, 1996.

\bibitem[Stoelinga(2002)]{stoelinga2002introduction}
Mari{\"e}lle Stoelinga.
\newblock An introduction to probabilistic automata.
\newblock \emph{Bulletin of the EATCS}, 78\penalty0 (176-198):\penalty0 2,
  2002.

\bibitem[Camacho and Alba(2013)]{camacho2013model}
Eduardo~F Camacho and Carlos~Bordons Alba.
\newblock \emph{Model predictive control}.
\newblock Springer science \& business media, 2013.

\bibitem[Bertsekas(2011)]{bertsekas2011dynamic}
Dimitri~P Bertsekas.
\newblock Dynamic programming and optimal control 3rd edition, volume ii.
\newblock \emph{Belmont, MA: Athena Scientific}, 2011.

\bibitem[Aguilar-Garcia et~al.(2015)Aguilar-Garcia, Fortes, Molina-Garc{\'\i}a,
  Calle-S{\'a}nchez, Alonso, Garrido, Fern{\'a}ndez-Dur{\'a}n, and
  Barco]{aguilar2015location}
Alejandro Aguilar-Garcia, S~Fortes, Mariano Molina-Garc{\'\i}a, Jaime
  Calle-S{\'a}nchez, Jos{\'e}~I Alonso, Aaron Garrido, Alfonso
  Fern{\'a}ndez-Dur{\'a}n, and Raquel Barco.
\newblock Location-aware self-organizing methods in femtocell networks.
\newblock \emph{Computer Networks}, 93:\penalty0 125--140, 2015.

\bibitem[Siddique et~al.(2015)Siddique, Tabassum, Hossain, and
  Kim]{siddique2015wireless}
Uzma Siddique, Hina Tabassum, Ekram Hossain, and Dong~In Kim.
\newblock Wireless backhauling of 5g small cells: Challenges and solution
  approaches.
\newblock \emph{IEEE Wireless Communications}, 22\penalty0 (5):\penalty0
  22--31, 2015.

\bibitem[Shakeri et~al.(2019)Shakeri, Al-Garadi, Badawy, Mohamed, Khattab,
  Al-Ali, Harras, and Guizani]{shakeri2019design}
Reza Shakeri, Mohammed~Ali Al-Garadi, Ahmed Badawy, Amr Mohamed, Tamer Khattab,
  Abdulla~Khalid Al-Ali, Khaled~A Harras, and Mohsen Guizani.
\newblock Design challenges of multi-uav systems in cyber-physical
  applications: A comprehensive survey and future directions.
\newblock \emph{IEEE Communications Surveys \& Tutorials}, 21\penalty0
  (4):\penalty0 3340--3385, 2019.

\bibitem[{Srivastava} and {Salapaka}(2020)]{9096570}
A.~{Srivastava} and S.~M. {Salapaka}.
\newblock Simultaneous facility location and path optimization in static and
  dynamic networks.
\newblock \emph{IEEE Transactions on Control of Network Systems}, pages 1--1,
  2020.

\bibitem[Mahajan et~al.(2009)Mahajan, Nimbhorkar, and
  Varadarajan]{mahajan2009planar}
Meena Mahajan, Prajakta Nimbhorkar, and Kasturi Varadarajan.
\newblock The planar k-means problem is np-hard.
\newblock In \emph{International Workshop on Algorithms and Computation}, pages
  274--285. Springer, 2009.

\bibitem[Srivastava and Salapaka(2021)]{9517030}
Amber Srivastava and Srinivasa~M. Salapaka.
\newblock Parameterized mdps and reinforcement learning problems--a maximum
  entropy principle-based framework.
\newblock \emph{IEEE Transactions on Cybernetics}, pages 1--13, 2021.
\newblock \doi{10.1109/TCYB.2021.3102510}.

\bibitem[Jaynes(2003)]{jaynes2003probability}
Edwin~T Jaynes.
\newblock \emph{Probability theory: The logic of science}.
\newblock Cambridge university press, 2003.

\bibitem[Rose(1991)]{rose1991deterministic}
Kenneth Rose.
\newblock \emph{Deterministic annealing, clustering, and optimization}.
\newblock PhD thesis, California Institute of Technology, 1991.

\bibitem[Chen et~al.(2005)Chen, Zhou, and Tang]{chen2005protein}
Luonan Chen, Tianshou Zhou, and Yun Tang.
\newblock Protein structure alignment by deterministic annealing.
\newblock \emph{Bioinformatics}, 21\penalty0 (1):\penalty0 51--62, 2005.

\bibitem[Xu et~al.(2014)Xu, Salapaka, and Beck]{xu2014aggregation}
Yunwen Xu, Srinivasa~M Salapaka, and Carolyn~L Beck.
\newblock Aggregation of graph models and markov chains by deterministic
  annealing.
\newblock \emph{IEEE Transactions on Automatic Control}, 59\penalty0
  (10):\penalty0 2807--2812, 2014.

\bibitem[Sharma et~al.(2012)Sharma, Salapaka, and Beck]{sharma2012entropy}
Puneet Sharma, Srinivasa~M Salapaka, and Carolyn~L Beck.
\newblock Entropy-based framework for dynamic coverage and clustering problems.
\newblock \emph{IEEE Transactions on Automatic Control}, 57\penalty0
  (1):\penalty0 135--150, 2012.

\bibitem[Xu et~al.(2013{\natexlab{a}})Xu, Salapaka, and Beck]{xu2013clustering}
Yunwen Xu, S~M Salapaka, and Carolyn~L Beck.
\newblock Clustering and coverage control for systems with acceleration-driven
  dynamics.
\newblock \emph{IEEE Transactions on Automatic Control}, 59\penalty0
  (5):\penalty0 1342--1347, 2013{\natexlab{a}}.

\bibitem[Srivastava and Salapaka(2020)]{srivastava2020simultaneous}
Amber Srivastava and Srinivasa~M Salapaka.
\newblock Simultaneous facility location and path optimization in static and
  dynamic networks.
\newblock \emph{IEEE Transactions on Control of Network Systems}, 2020.

\bibitem[Xu et~al.(2013{\natexlab{b}})Xu, Salapaka, and Beck]{xu2013distance}
Yunwen Xu, Srinivasa~M Salapaka, and Carolyn~L Beck.
\newblock A distance metric between directed weighted graphs.
\newblock In \emph{52nd IEEE Conference on Decision and Control}, pages
  6359--6364. IEEE, 2013{\natexlab{b}}.

\bibitem[Biondi et~al.(2014)Biondi, Legay, Nielsen, and
  Wkasowski]{biondi2014maximizing}
Fabrizio Biondi, Axel Legay, Bo~Friis Nielsen, and Andrzej Wkasowski.
\newblock Maximizing entropy over markov processes.
\newblock \emph{Journal of Logical and Algebraic Methods in Programming},
  83\penalty0 (5-6):\penalty0 384--399, 2014.

\bibitem[Savas et~al.(2018)Savas, Ornik, Cubuktepe, and
  Topcu]{savas2018entropy}
Yagiz Savas, Melkior Ornik, Murat Cubuktepe, and Ufuk Topcu.
\newblock Entropy maximization for constrained markov decision processes.
\newblock In \emph{2018 56th Annual Allerton Conference on Communication,
  Control, and Computing}, pages 911--918. IEEE, 2018.

\bibitem[Srivastava and Salapaka(2022)]{srivastava2021paraSDM}
Amber Srivastava and Srinivasa~M Salapaka.
\newblock Parameterized mdps and reinforcement learning problems - a maximum
  entropy principle based framework.
\newblock \emph{arXiv preprint arXiv:2006.09646}, 2022.

\bibitem[Sontag(1983)]{sontag1983lyapunov}
Eduardo~D Sontag.
\newblock A lyapunov-like characterization of asymptotic controllability.
\newblock \emph{SIAM Journal on Control and Optimization}, 21\penalty0
  (3):\penalty0 462--471, 1983.

\bibitem[Sepulchre et~al.(2012)Sepulchre, Jankovic, and
  Kokotovic]{sepulchre2012constructive}
Rodolphe Sepulchre, Mrdjan Jankovic, and Petar~V Kokotovic.
\newblock \emph{Constructive nonlinear control}.
\newblock Springer Science \& Business Media, 2012.

\bibitem[Royden and Fitzpatrick(1988)]{royden1988real}
Halsey~Lawrence Royden and Patrick Fitzpatrick.
\newblock \emph{Real analysis}, volume~32.
\newblock Macmillan New York, 1988.

\end{thebibliography}

\end{document}